\documentclass{article}
\usepackage{amsmath}
\usepackage{amsfonts}
\usepackage{amssymb}	
\usepackage{color}
\newtheorem{theo}{Theorem}[section]

\newtheorem{rem}[theo]{Remark}
\newtheorem{prop}[theo]{Proposition}

\newcommand{\mysection}[1]{\section{#1} \setcounter{equation}{0}}
\newcommand{\proof}{{\sc Proof.} \quad}

\newcommand{\be}{\begin{equation} \label}
\newcommand{\ee}{\end{equation}}
\newcommand{\bes}{\begin{equation} \begin{array}{c} \label}
\newcommand{\ees}{\end{array} \end{equation}}
\newcommand{\bea}{\begin{eqnarray}\label}
\newcommand{\eea}{\end{eqnarray}}
\newcommand{\beas}{\begin{eqnarray} \begin{array}{rcl} \label}
\newcommand{\eeas}{\end{array} \end{eqnarray}}
\newcommand{\bas}{\begin{eqnarray*}}\newcommand{\eas}{\end{eqnarray*}}
\newcommand{\bass}{\begin{eqnarray*} \begin{array}{rcl}}
\newcommand{\eass}{\end{array} \end{eqnarray*}}
\newcommand{\basss}{\begin{eqnarray*} \begin{array}{c}}
\newcommand{\easss}{\end{array} \end{eqnarray*}}
\newcommand{\qed}{{}\hfill $\square$ \\}
\newcommand{\bit}{\begin{itemize}}
\newcommand{\eit}{\end{itemize}}

\begin{document}
\title{Unboundedness of solutions to a supercritical quasilinear parabolic-parabolic Keller-Segel system in dimension $1$ and to a critical one in higher dimensions}
\author{
Tomasz Cie\'{s}lak \\
{\small Institute of Mathematics, Polish Academy of Sciences, \'Sniadeckich 8, 00-956 Warsaw, Poland}\\
{\small E-Mail: T.Cieslak@impan.pl}}

\maketitle

\begin{abstract}
Unboundedness of solutions is shown to occur in a one-dimensional quasilinear parabolic-parabolic chemotaxis system for any initial mass. Our result is also independent of the relation between the speeds of the diffusion of cells and chemoattractant. The proof is achieved by contradiction. It uses a virial type inequality together with a boundedness from below of the Liapunov functional associated to the system. Moreover, the same method is applied also in higher dimensions and an unboundedness result is shown in the case of critical quasilinear fully parabolic Keller-Segel system
for mass large enough in dimensions $n\geq 3$.  

\end{abstract}

\mysection{Introduction}\label{section1}
This work deals with nonnegative solution couples $(u,v)$ of the parabolic-parabolic 
Keller-Segel system
\be{0}
	\left\{ \begin{array}{ll}
	u_t= (a(u)u_x)_x- (uv_x )_x, & \; x\in (0,1), \ t>0, \\[2mm]
	\tau v_t=v_{xx}-v+u, & \; x\in (0,1), \ t>0, \\[2mm]
	u_x(0)=u_x(1)=v_x(0)=v_x(1)=0, &\; t>0, \\[2mm]
	u(x,0)=u_0(x), \quad v(x,0)=v_0(x), & \; x\in (0,1),
	\end{array} \right.
\ee
where $\tau>0$
and the initial data are supposed to satisfy $u_0, v_0\in W^{1,p}(0,1)$, $p>1$ such that
$u_0 \geq 0$ and $v_0 \geq 0$ in $(0,1)$.
Moreover, we assume that $a\in C^2([0,\infty))$ is a positive function.

Systems of this kind were introduced in \cite{KS:1} to describe the motion of cells which are diffusing and moving towards the gradient of a substance called chemoattractant, the latter being produced by cells themselves. The original model consisted of several equations including more unknowns. It was simplified in \cite{Nan} where the form of the equation similar to \eqref{0} was introduced. The essentiality of the nonlinear diffusion was emphasized in \cite{hil:volume} where the so-called volume-filling effect was studied. The important issue in studies concerning \eqref{0} was a chemotactic collapse of cells interpreted as a finite-time blowup of the component $u$. In \cite{child},
based on the numerical simulations, the authors suggested that in a two-dimensional (most popular and biologically most relevant) semilinear version, there is a threshold of mass such that a solution starting from the initial data with a mass smaller than a threshold exists globally, thus avoiding the chemotactic collapse, while initial conditions having mass larger than a threshold lead to a chemotactic collapse. We emphasize at this point that the rigorous studies that followed were related to further simplifications of the Keller-Segel model. Actually, until recently, almost all the results concerning finite-time blowup were achieved for a parabolic-elliptic simplification ($\tau=0$ in \eqref{0}). First, in \cite{jl:expl} a conjecture from \cite{child} was shown in the two-dimensional semilinear parabolic-elliptic simplification. In \cite{Na95} the exact threshold in the case of radially symmetric data was calculated as $8\pi$ in the parabolic-elliptic 2d case. Next, the critical mass phenomenon was also proved in the case of a two-dimensional problem for solutions without symmetries (see \cite{Nagai}, \cite{NSY97}). Moreover, it has been shown that in higher dimensions a finite-time blowup of solutions to the semilinear version of \eqref{0} can occur independently of the initial mass provided that the initial data are concentrated enough \cite{Na95}.  Finally, during last years, the quasilinear system was studied. Critical nonlinearities (and the lack of it in dimension $1$) and critical mass phenomenons were found in all the dimensions, see \cite{CW08},\cite{clGlob_vs_blow},\cite{CRAS},\cite{Nasri}.  

However, all those results are available only for a parabolic-elliptic simplification of \eqref{0}. In the case of the original fully parabolic version the investigation of chemotactic collapse turned out to be a much more challenging issue. Until 2010 the only available result stating the occurrence of finite-time blowup was the one in \cite{HV:1}, where examples of some particular solutions to the semilinear version of \eqref{0} in dimension $n=2$ blowing up in a finite-time are shown. Moreover, see \cite{Ho02}, \cite{ss1}, it was known that solutions to a fully parabolic problem can be unbounded. However it was not stated if blowup occurs in finite or infinite time. The same approach was used in \cite{win_mmas} to show unboundedness of solutions to a fully parabolic Keller-Segel system in the case of supercritical quasilinear systems in dimensions $n\geq 2$. Next, in \cite{clKS} the collapse of solutions to the one-dimensional Keller-Segel system with appropriately weak diffusion of cells and sufficiently fast diffusion of chemoattractant is shown provided initial mass is large enough. The breakthrough has been made recently in \cite{win_bu}. Introducing a new method M. Winkler
shows there that in dimensions $n\geq 3$ generic solutions to the semilinear version of \eqref{0} blow up in finite time independently of the size of initial mass. His method was generalized to the quasilinear case in \cite{ciesti}, \cite{ciesti2}, where the finite-time blowup of solutions to a fully parabolic supercritical Keller-Segel system was shown for any initial mass in any dimension $n\geq 2$. The problem for critical nonlinearities is still open.

This way, the occurence of finite-time singularities of solutions to a fully parabolic quasilinear Keller-Segel system in higher dimensions is already known for supercritical diffusions independently of the size of initial mass. However, in the one-dimensional case, the only known result is the one from \cite{clKS}, where the finite-time blowup is shown only for large mass initial data and moreover in a neighbourhood of a parabolic-elliptic situation, meaning $\tau$ being a very small number. The purpose of the present paper is to introduce a result concerning the supercritical case of \eqref{0} with any $\tau>0$ and an arbitrary mass. We are going to give a result stating the unboundedness of solutions. It does not imply a finite-time blowup. We do not exclude a situation that a solution becomes unbounded when time goes to infinity. We know that this kind of situations might happen (see \cite[Theorem 1.6]{ciesti} or \cite[Theorem 1.4]{ciesti2}) in chemotaxis. Our result is similar to the higher-dimensional results in \cite{ss1}, \cite{Ho02} or \cite{win_mmas}, however the method of a proof which we follow is totally different. The method of \cite{ss1}, \cite{Ho02}, \cite{win_mmas} is not applicable in the one-dimensional case. Indeed, unboundedness of solutions of a higher dimensional Keller-Segel system in the above-mentioned articles is a consequence of a contradiction between boundedness of solutions which immediately leads to convergence with time (at least on subsequences) to steady states and the fact that one can construct  initial conditions such that the value of the Liapunov functional on them is smaller than the minimum of this Liapunov functional over the set of steady states. This approach fails in the one-dimensional case due to the fact that the 1d Liapunov functional is bounded from below by its minimum over the steady states. We develop an alternative way of proving unboundedness of solutions in that case.  We use an extension of a moment method (see \cite{BiN94}, \cite{Na95} or \cite{clSP}). 

The same method can be applied also in higher dimensions. However, in supercritical cases, unboundedness results, as well as much stronger finite-time blowup results, are already known, see \cite{Ho02}, \cite{ss1}, \cite{win_mmas}, \cite{win_bu}, \cite{ciesti}.  But so far, there were no results (neither those concerning finite nor infinite-time blowup) in the critical case in dimensions $n\geq 3$. We fill in this gap in the last section.

\mysection{Preliminaries}\label{section2} 

First, notice that thanks to the application of Amann's theory (\cite{amann}) we have a local well-posedness result for the regular nonnegative solutions, once we start from the regular nonnegative initial data, see (\cite{clKS}). Integrating by parts we see that those solutions also satisfy 
\begin{equation}\label{mass}
\int_0^1u(t,x)dx=\int_0^1u_0(x)dx, \int_0^1v(t,x)dx=\int_0^1u_0(x)dx+e^{-\frac{t}{\tau}}\left(\int_0^1v_0(x)dx-\int_0^1u_0(x)dx\right)
\end{equation}
for any $t>0$. Moreover, if the solution ceases to exist for times larger than $t_0$, then immediately $L^\infty$ norm of $u$ blows up at a time $t_0$.

Next we recall the existence and properties of a Liapunov functional associated to \eqref{0}.
\begin{equation}\label{b6}
L(u(t),v(t)) + \tau\int_0^t \|\partial_t v(s)\|_2^2\ ds+ \int_0^t\int_0^1\frac{(a(u)u_x-uv_x)^2}{u}=L(u_0,v_0) \;\;\mbox{ for }\;\; t>0\,,
\end{equation}
where
\begin{equation}\label{b7}
L(u,v):= \int_0^1 \left( \Phi(u) - u v + \frac{1}{2}\ |\partial_x v|^2 + \frac{1}{2}\ |v|^2 \right)\ dx\,,
\end{equation}
where $\Phi\in \mathcal{C}^2((0,\infty))$ is defined in a following way
\begin{equation}\label{b2}
\Phi(1)=\Phi'(1):=0 \;\;\mbox{ and }\;\; \Phi''(r) := \frac{a(r)}{r} \;\;\mbox{ for }\;\; r>0\,.
\end{equation}
We take an advantage of the one-dimensional setting and state the boundedness from below of the Liapunov functional, compare \cite[Lemma 5]{clKS}
\begin{equation}\label{b8}
L(u(t),v(t)) \ge - C \;\;\mbox{ for }\;\; t>0\,.
\end{equation}
The main theorem of this part reads
\begin{theo}\label{Tw.1}
Assume $a\in L^1(0,\infty)$. Then for any $\tau>0$ and any initial mass $m>0$ there exists an initial data $(u_0(x),v_0(x))$ such that $\int_0^1u_0(x)dx=m$ and a component $u(t)$ of a solution to \eqref{0} starting from $(u_0, v_0)$ is unbounded in an $L^\infty$ norm.  
\end{theo}
\begin{rem}\label{Uwaga1}
Our theorem does not say whether the blowup occurs in finite or infinite time.
\end{rem}
In the sequel we will need the following observation which is a consequence of classical regularity estimates for parabolic equations, however in our case we can infer it easily, keeping the proof self-consistent.
\begin{prop}\label{Prop1}
Assume $\left\|u(t,\cdot)\right\|_{L^\infty(0,1)}<C$ for $t>0$, then there exists a constant $c_2>0$ such that $\left\|v_t(t,\cdot)\right\|_{2}<c_2$ for $t>0$.
\end{prop} 
Proof.
First we differentiate the second equation in \eqref{0} in $t$ and next multiply it by $v_t$ and integrate. We arrive at
\[
\frac{\tau}{2}\frac{d}{dt}\int_0^1v_t^2dx+\int_0^1(v_{tx})^2dx+\int_0^1v_t^2dx=\int_0^1u_tv_tdx.
\] 
By the first equation in \eqref{0} we see that
\[
\int_0^1u_tv_tdx=\int_0^1(a(u)u_x-uv_x)v_{tx}dx\leq \frac{1}{2}\int_0^1(a(u)u_x-uv_x)^2dx+ \frac{1}{2}\int_0^1(v_{tx})^2dx.
\] 
Hence
\begin{equation}\label{eins}
\frac{\tau}{2}\frac{d}{dt}\int_0^1v_t^2dx+\frac{1}{2}\int_0^1(v_{tx})^2dx+\int_0^1v_t^2dx\leq \frac{\left\|u\right\|_{L^\infty(0,1)}}{2}\int_0^1\frac{(a(u)u_x-uv_x)^2}{u}dx.
\end{equation}
Next we observe that due to (\ref{b6}) and (\ref{b8}) we have
\[
\frac{1}{2}\int_0^t\int_0^1\frac{(a(u)u_x-uv_x)^2}{u}dxdt<\infty\;\;\mbox{for any}\;t>0. 
\]
Thus, integrating (\ref{eins}) with respect to time, keeping in mind that $\left\|u(t,\cdot)\right\|_{L^\infty(0,1)}$ is assumed to be bounded, we arrive at the desired estimate of $v_t$.
\qed

\mysection{A proof of the main result}\label{section3}

This section is devoted to the proof of Theorem \ref{Tw.1}. We are going to use a method of a generalized second moment as introduced in \cite{clKS}. However,  our aim is weaker than in a mentioned paper. Thanks to the estimate on $v_t$ coming from Proposition \ref{Prop1}, by applying some ideas from \cite{clSP}, we arrive at the estimate which in turn allows us to infer unboundedness result for arbitrary initial mass and arbitrary $\tau>0$. This result is of special interest due to the fact that among supercritical Keller-Segel systems, the one-dimensional one was by now the only case where not only finite-time blowup, but even unboundedness of solutions was not known. In dimensions $n\geq 2$ it is already known in supercritical cases that for any mass there exist solutions blowing up in finite time, see \cite{ciesti}, \cite{ciesti2}. However, the method used in those cases, being an extension of the method of Winkler, \cite{win_bu}, does not seem to be applicable in the one-dimensional case due to the bound of the Liapunov functional from below, see (\ref{b8}). For the same reason the usual method of showing unboundedness of solutions to the higher-dimensional Keller-Segel system does not work in the 1d setting. This method, see \cite{Ho02}, \cite{ss1}, \cite{win_mmas}, requires unboundedness af a Liapunov functional in order to enable the user to choose initial data such that the value of the Liapunov functional on them is smaller than the infimum of the values of the Liapunov functional over the steady states. In a one-dimensional setting the minimal value of the Liapunov functional seems to be attained at the steady states, excluding this way the above method.

Our proof is of a reduction to absurd nature. 

\vspace{0.3cm}
\noindent Proof.
Let $q>2$ and define the generalized moment in a following way
\[
M_q (t) : =\frac{1}{q} \int_0^1 |U (t,x)|^q dx,
\] 
where 
\[
U (t, x) : = \int_0^x u(t, z) dz\;\;\mbox{for}\;  x \in [0, 1]\;\mbox{and}\;
V (t, x) := \int_0^x v(t, z) dz \;\mbox{for}\; x \in [0, 1].
\]
Assume $m:=\int_0^1u_0(x)dx>0$, 
\begin{equation}\label{v_0}
0<\int_0^1v_0(x)dx-m<\frac{m}{2(q+1)}.
\end{equation}
Moreover assume $\sup_{t>0}\left\|u(t,\cdot)\right\|_{L^\infty(0,1)}<\infty.$
By \eqref{0} we have
\begin{equation}\label{3.3.3}
\frac{d}{dt} M_q(t) =  \int_0^1 U^{q-1}(A(u))_x -  \int_0^1 U^{q-1}  u(\tau V_t-U+V),
\end{equation}
where $A(u):=-\int_u^\infty a(s)ds$.
Integrating by parts we have
\begin{multline}\label{3.3.4}
\frac{d}{dt}M_q(t) = - (q-1) \int_0^1 U^{q-2} u A(u) +(U^{q-1} A(u))\Big|_0^1 +   \int_0^1 U^{q-1}  u(U-V-\tau V_t)  =\\ (q-1) \int_0^1 U^{q-2} (-u A(u)) + m^{q-1} A(u (t,1))  + \frac{1}{q+1} \int_0^1 (U^{q+1})_x - \frac{1}{q} \int_0^1 (U^{q})_x (V+\tau V_t).
\end{multline}
Due to (\ref{mass}) and (\ref{v_0}) we have
\[
-\frac{1}{q}\int_0^1(U^q)_xV=\frac{1}{q}\int_0^1U^qv-\frac{1}{q}(m^qV(t,1))\leq \frac{1}{q}\int_0^1U^qv-\frac{1}{q}m^{q+1},
\]
which in turn, together with (\ref{3.3.3}) and an observation $A(u (t,1)) \le 0$ yields
\begin{equation}\label{3-3}
\frac{d}{dt}M_q(t)\leq (q-1) \int_0^1 U^{q-2} (-u A(u))+\left\|v\right\|_{L^\infty(0,1)}M_q(t)-\frac{1}{q(q+1)}m^{q+1}+\frac{\tau }{q} \left[ \left( \int_0^1 U^q v_t \right) -  m^q V_t  (1) \right].
\end{equation}
Since 
\[
-V_t(1)=\frac{1}{\tau}e^{-\frac{t}{\tau}}\left(\int_0^1v_0dx-m\right),
\] 
in view of (\ref{v_0}) we have
\[
\frac{d}{dt}M_q(t)\leq (q-1) \int_0^1 U^{q-2} (-u A(u))+\left\|v\right\|_{L^\infty(0,1)}M_q(t)-\frac{m^{q+1}}{q(q+1)}+\frac{\tau }{q}\left( \int_0^1 U^q v_t \right)
+\frac{m^q}{q}\left(\int_0^1v_0dx-m\right)
\]
\[
\leq(q-1) \int_0^1 U^{q-2} (-u A(u))+\left\|v\right\|_{L^\infty(0,1)}M_q(t)-\frac{1}{2q(q+1)}m^{q+1}+\frac{\tau }{q}\left( \int_0^1 U^q v_t \right).
\]
Next, by H\"{o}lder's inequality 
\begin{equation}\label{tc}
\frac{d}{dt}M_q(t)\leq (q-1) \int_0^1 U^{q-2} (-u A(u))+\left\|v\right\|_{L^\infty(0,1)}M_q(t)-\frac{1}{2q(q+1)}m^{q+1}+\frac{\tau m^{q/2}}{q^{1/2}}M_q^{\frac{1}{2}}\left\|v_t\right\|_2.
\end{equation}
We are now in a position to introduce a concave function $B$ such that
\begin{eqnarray}
0 \le -r A(r) & \le & B(r) \,, \qquad r\in (0,\infty)\,, \label{bu1} \\
\lim_{r\to \infty} \frac{B(r)}{r} & = & 0\,. \label{bu2}
\end{eqnarray}
This function was introduced in \cite{clSP} and in \cite[Lemma 3.1]{clGlob_vs_blow} it was shown that such a function can be associated to any $a$ satisfying the assumptions of Theorem \ref{Tw.1}. Its advantage is that we can proceed for any $m>0$. 
The following estimates, using in a tricky way Jensen's inequality for several probabilistic measures, 
are taken from \cite[Theorem 10]{clSP}, we provide them for the sake of completeness.
\[
\int_0^1 U^{q-2} B(u)dx =\left( \int_0^1 B(u) dx\right) \int_0^1 (U^{q})^\frac{q-2}{q} \frac{B(u) dx}{\int_0^1 B(u)dx}
\] 
\[
\le \left( \int_0^1 B(u) dx\right)  \left(\int_0^1 U^{q} \frac{B(u) dx}{\int_0^1 B(u)dx} \right)^\frac{q-2}{q} 
\] 
\[
= \left( \int_0^1 B(u) dx\right)^\frac{2}{q} (qM_q)^\frac{q-2}{q}   \left(\int_0^1 B(u) \frac{U^q dx}{q M_q} \right)^\frac{q-2}{q}
\]
\begin{equation}\label{3.3.5} 
\le 
 (qM_q)^\frac{q-2}{q}  \left[ B \left( \int_0^1 u \right) \right]^\frac{2}{q}  \left[ B \left(  \int_0^1 \frac{u U^q dx}{q M_q} \right) \right]^\frac{q-2}{q}\leq \left(qM_q\right)^{\frac{q-2}{2}}[B(m)]^{2/q}\left[B\left(\frac{m^{q+1}}{q(q+1)M_q}\right)\right]^{\frac{q-2}{2}}.
\end{equation}
The Jensen inequality above was applied to $\frac{B(u) dx}{\int B(u)}$, $dx$ and $\frac{U^q dx}{q M_q}$. 

In view of (\ref{3.3.5}), Proposition \ref{Prop1} and the bound on $\left\|v\right\|_{L^\infty(0,1)}$, (\ref{tc}) reads
\begin{equation}\label{tc1}
\frac{d}{dt} M_q(t)\leq (q-1) (qM_q)^\frac{q-2}{q} [B(m)]^\frac{2}{q}  \left[ B \left(  \int_0^1 \frac{u U^q dx}{q M_q} \right) \right]^\frac{q-2}{q}+c_1M_q(t)-\frac{m^{q+1}}{2q(q+1)}+\frac{c_2\tau m^{q/2}}{q^{1/2}}M_q^{\frac{1}{2}}.
\end{equation}
Next (\ref{tc1}) reads
\begin{equation}\label{tc3}
\frac{d}{dt} M_q(t)\leq \Lambda_m(M_q(t)),
\end{equation}
where 
\[
\Lambda(r):=c_1r+(q-1)[B(m)]^{2/q}\left(\frac{m^{q+1}}{q+1}\right)^{\frac{q-2}{2}}\left[\beta\left(\frac{m^{q+1}}{q(q+1)r}\right)\right]^{\frac{q-2}{2}}+\frac{c_2\tau m^{q/2}}{q^{1/2}}r^{\frac{1}{2}}-\frac{m^{q+1}}{2q(q+1)}
\]
for $r>0$ and $\beta(r)=\frac{B(r)}{r}$.

We notice that according to (\ref{bu2}) for any $m>0$ $\Lambda_m(0)<0$. Due to continuity of $\Lambda$, there exists $\theta>0$ such that $\Lambda(r)<0$ for $r\in[0,\theta]$. Thus, being allowed to choose initial data such that $M_q(0)<\theta$ we arrive at the claim that $M_q$ must touch $0$ at a finite time. But this is contradictory to the nonnegativity of $u$ and a mass conservation (\ref{mass}), what means that a solution ceases to exist at latest at the time of an extinction of $M_q$. This in turn means that $u$ is not bounded in $L^\infty$, a contradiction. 

\qed 

\mysection{An application of the method in higher dimensions} 
This section is devoted to proving unboundedness of radially symmetric solutions to the fully parabolic quasilinear critical Keller-Segel system in dimensions $n\geq 3$. We achieve this issue extending the method of previous section to higher dimensions. Again, our result does not say whether blowup of solutions occurs in finite or infinite time.

As is well known (see for example \cite{CW08}) radial symmetry is preserved by the Keller-Segel model in a ball and equations of evolution in that case read
\be{OR}
  \left\{
  \begin{array}{l}
    u_t=\frac{1}{r^{n-1}} (r^{n-1} a(u) u_r)_r
    - \frac{1}{r^{n-1}} (r^{n-1} uv_r)_r, \qquad r \in (0,1), \ t>0, \\[2mm]
    v_t = \frac{1}{r^{n-1}} (r^{n-1} v_r)_r + u - v, \qquad r \in (0,1), \ t>0, \\[2mm]
    u_r=v_r=0, \qquad r\in\{0,1\}, \ t>0, \\[2mm]
    u(r,0)=u_0(r), v(r,0)=v_0(r)\qquad r \in (0,1).
  \end{array} \right.
\ee
Like in the case of (\ref{0}) $u,v$ denote the density of cells and chemoattractant. Local-in-time solutions exist and whenever they cease to exist for all times we know that $\left\|u(t,\cdot)\right\|_{L^\infty(B(0,1))}$ goes to infinity for some $t<\infty$.  Moreover both $u_0, v_0$ are assumed to be nonnegative, hence for $t>0$ both $u(t,\cdot), v(t,\cdot)$ are positive due to a comparison principle. The function $a$ is $C^2$ regular and positive, it stands for nonlinear diffusion. The radius of a ball is $1$, and the radial variable is denoted by $r: 0<r<1$. Next,
\be{mass1}
	\left\{
	\begin{array}{l}
		\int_{B(1,0)}u(t,x)dx=\int_{B(1,0)}u_0(x)dx,\\[2mm] 					
		\int_{B(1,0)}v(t,x)dx=\int_{B(1,0)}u_0(x)dx+e^{-t}\left(\int_{B(1,0)}v_0(x)dx-\int_{B(1,0)}u_0(x)dx\right)\\[2mm]
	\end{array}\right.
\ee  
hold, where $B(0,1)$ is a ball centred at $0$. 

Like in an one-dimensional case, we have
\begin{equation}\label{b16}
L(u(t),v(t)) + \int_0^t \|\partial_t v(s)\|_2^2\ ds+ \int_0^t\int_{B(0,1)}\frac{|a(u)\nabla u-u\nabla v|^2}{u}=L(u_0,v_0) \;\;\mbox{ for }\;\; t>0\,,
\end{equation}
with $L$ defined exactly as in (\ref{b7}) with integration over the interval replaced by the one over a ball.
Moreover, we have a simple proposition
\begin{prop}\label{Prop11}
Assume $\left\|u(t,\cdot)\right\|_{L^\infty(B(0,1))}<C$ for $t>0$, then there exists $C_1>0$ such that $L\geq -C_1$.
\end{prop}
This is implied by the fact that $\int_{B(0,1)}uv$ is bounded independently of time. Next, analogously to Proposition
\ref{Prop1} we have
\begin{prop}\label{Prop12}
Assume $\left\|u(t,\cdot)\right\|_{L^\infty(B(0,1))}<C$ for $t>0$, then $\left\|v_t(t,\cdot)\right\|_{2}<C_1$ for $t>0$.
\end{prop}
We are now in a position to formulate the main theorem of this section.
\begin{theo}\label{Tw.11}
Let $n\geq 3$. Assume $a(u)=(1+u)^{1-\frac{2}{n}}$. Then there exists such value, say $m_*$, that for any $M>m_*$ there exists $(u_0, v_0)$ with mass $\int_{B(0,1)}u_0(x)dx=M|B(0,1)|$, such that a component $u(t)$ of a solution of \eqref{OR} starting from $u_0$ is unbounded in $L^\infty$ norm.  
\end{theo}
\proof Like in section \ref{section3} first we assume that $u$ is bounded in $L^\infty$. In what follows we show that this assumption leads to an absurd. Let us follow the ideas which appeared in \cite{CRAS}. We introduce
\[
U(t, r) := \frac{1}{n|B(0, 1)|}\int_{B(0,r)}u(t,x) dx\;\;\mbox{and}\;\; M_2(t) := \frac{1}{2}\int_0^1\left(\frac{M}{n}-U(t,r)\right)^{2}r^{n-1}dr
\]
for $(t,r)\in [0,T_{max})\times[0,1]$. Next, we notice that
\begin{equation}\label{napisy}
U(t,0)=\frac{M}{n}-U(t,1)=0,
\end{equation}
Moreover for $V$ defined analogously to $U$ we have
\begin{equation}\label{napisy1}
V(t,1)=\frac{M}{n}+e^{-t}\left(V(0,1)-\frac{M}{n}\right).
\end{equation}
Integrating (\ref{OR}) we get
\[
U_t=r^{n-1}A(u)_r-r^{n-1}uv_r,\;\;V_t=r^{n-1}v_r-V+U,
\]
where $A$ is a primitive of $a$ such that $A(0)=0$. Hence
\[
\frac{dM_2}{dt}=\int_0^1\left(\frac{M}{n}-U\right)\left[-r^{2(n-1)}(A(u))_r-\frac{M}{n}U_r+\left(\frac{M}{n}-U\right)U_r+VU_r+V_tU_r\right]dr.
\]
Consequently, 
\[
\frac{dM_2}{dt}=\int_0^1\left[2(n-1)r^{2n-3}\left(\frac{M}{n}-U\right)-r^{2(n-1)}U_r\right]A(u)dr-\frac{1}{6}\left(\frac{M}{n}\right)^3+
\]
\[
-\int_0^1\left[\frac{1}{2}\left(\frac{M}{n}-U\right)^2\right]_r(V+V_t)dr.
\]
Thus, due to the fact that integration by parts gives
\[
-\int_0^1\left[\frac{1}{2}\left(\frac{M}{n}-U\right)^2\right]_r(V+V_t)dr=\int_0^1\frac{1}{2}\left(\frac{M}{n}-U\right)^2r^{n-1}(v_t+v)dr,
\]
we infer
\begin{equation}\label{oszacowanie}
\frac{dM_2}{dt}\leq \int_0^1\left[2(n-1)r^{2n-3}\left(\frac{M}{n}-U\right)-r^{2(n-1)}U_r\right]A(u)dr-\frac{1}{6}\left(\frac{M}{n}\right)^3
+\left\|v\right\|_{L^\infty(B(0,1))}M_2(t)+
\end{equation}
\[
\frac{M}{n}(M_2(t))^{\frac{1}{2}}\left(\int_0^1v_t^2r^{n-1}dr\right)^{\frac{1}{2}}
\]
\[
=\int_0^1\left[2(n-1)r^{2n-3}\left(\frac{M}{n}-U\right)-r^{2(n-1)}U_r\right]A(u)dr-\frac{1}{6}\left(\frac{M}{n}\right)^3
\]
\[
+\left\|v\right\|_{L^\infty(B(0,1))}M_2(t)+\frac{M}{n^{3/2}|B(0,1)|^{1/2}}(M_2(t))^{\frac{1}{2}}\left\|v_t\right\|_{L^2(B(0,1))}^{\frac{1}{2}}
\]
Keeping in mind that $a(u)=(1+u)^{1-\frac{2}{n}}$, we see that
\[
A(u)\leq \frac{1}{2-\frac{2}{n}}u^{2-\frac{2}{n}}+u.
\]
Hence we are now in a position to use computations from \cite{CRAS} in order to estimate a term 
\begin{equation}\label{term}
\int_0^1\left[2(n-1)r^{2n-3}\left(\frac{M}{n}-U\right)-r^{2(n-1)}U_r\right]A(u)dr.
\end{equation}
For the sake of completeness we provide them with all details. First we notice that
\[
\left[2(n-1)\left(\frac{M}{n}-U\right)-r^nu\right]A(u)\leq \max\{2(n-1)\left(\frac{M}{n}-U\right)-r^nu,0\}\left(\frac{1}{2-\frac{2}{n}}u^{2-\frac{2}{n}}+u\right).
\] 
Next, 
\[
\max\{2(n-1)\left(\frac{M}{n}-U\right)-r^nu,0\}\left(\frac{1}{2-\frac{2}{n}}u^{2-\frac{2}{n}}+u\right)\leq
\]
\[
\max\{2(n-1)\left(\frac{M}{n}-U\right)-r^nu,0\}\left(\frac{1}{2-\frac{2}{n}}\left[\frac{2(n-1)}{r^n}\left(\frac{M}{n}-U\right)\right]^{1-\frac{2}{n}}+1\right)u\leq
\]
\[
\frac{1}{2-\frac{2}{n}}\left[2(n-1)\left(\frac{M}{n}-U\right)\right]^{2-\frac{2}{n}}r^{-n+2}u+2(n-1)\left(\frac{M}{n}-U\right)u,
\]
in turn (\ref{term}) is estimated by
\[
\int_0^1\frac{1}{2-\frac{2}{n}}\left[2(n-1)\left(\frac{M}{n}-U\right)\right]^{2-\frac{2}{n}}U_rdr+ 2(n-1)\int_0^1r^{n-2}\left(\frac{M}{n}-U\right)U_rdr=
\]
\begin{equation}\label{term1}
\frac{(2(n-1))^{2-\frac{2}{n}}}{(2-\frac{2}{n})(3-\frac{2}{n})}\left(\frac{M}{n}\right)^{3-\frac{2}{n}}+ 2(n-1)\int_0^1r^{n-2}\left(\frac{M}{n}-U\right)U_rdr.
\end{equation}
A term $2(n-1)\int_0^1r^{n-2}\left(\frac{M}{n}-U\right)U_rdr$ can be estimated further using the fact that the function
$r\rightarrow r^{1-\frac{2}{n}}$ is concave. We use the Jensen inequality applied to a measure $\left(\frac{M}{n}-U\right)U_rdr$, integrate by parts and obtain
\begin{equation}\label{term2}
2(n-1)\int_0^1r^{n-2}\left(\frac{M}{n}-U\right)U_rdr\leq (n-1)(2n)^{1-\frac{2}{n}}\left(\frac{M}{n}\right)^{\frac{4}{n}}M_2(t)^{1-\frac{2}{n}}.
\end{equation}
We conclude that in view of (\ref{term1}) and (\ref{term2}), (\ref{term}) is estimated by
\[
\frac{(2(n-1))^{2-\frac{2}{n}}}{(2-\frac{2}{n})(3-\frac{2}{n})}\left(\frac{M}{n}\right)^{3-\frac{2}{n}}+(n-1)(2n)^{1-\frac{2}{n}}\left(\frac{M}{n}\right)^{\frac{4}{n}}M_2(t)^{1-\frac{2}{n}},
\]
hence (\ref{oszacowanie}) yields 
\begin{equation}\label{istotne}
\frac{dM_2(t)}{dt}\leq \frac{(2(n-1))^{2-\frac{2}{n}}}{(2-\frac{2}{n})(3-\frac{2}{n})}\left(\frac{M}{n}\right)^{3-\frac{2}{n}}-\frac{1}{6}\left(\frac{M}{n}\right)^3+(n-1)(2n)^{1-\frac{2}{n}}\left(\frac{M}{n}\right)^{\frac{4}{n}}M_2(t)^{1-\frac{2}{n}}+
\end{equation}
\[
\left\|v\right\|_{L^\infty(B(0,1))}M_2(t)+\frac{M}{n^{3/2}|B(0,1)|^{1/2}}(M_2(t))^{\frac{1}{2}}\left\|v_t\right\|_{L^2(B(0,1))}^{\frac{1}{2}}.
\]
Now, due to an assumption that $L^\infty$ norm of $u$ is bounded and a local existence result, we infer that solution exists for any time $t>0$. However, as we see in a moment, for initial mass greater than $m_*$ we can find such initial data $(u_0,v_0)$, that this is in contradiction with (\ref{istotne}). Indeed, according to Proposition \ref{Prop12}  
$\left\|v_t\right\|_{L^2(B(0,1))}$ is bounded, moreover $v$ is also bounded in $L^\infty$ norm due to the standard parabolic regularity estimates related to the second equation in \eqref{OR}. Hence for $M$ large enough that 
\[
\frac{(2(n-1))^{2-\frac{2}{n}}}{(2-\frac{2}{n})(3-\frac{2}{n})}\left(\frac{M}{n}\right)^{3-\frac{2}{n}}-\frac{1}{6}\left(\frac{M}{n}\right)^3<0
\]
and choosing initial data $u_0$ concentrated enough that $M_2(u_0)$ is close to $0$, we see that the right-hand side of 
(\ref{istotne}) is negative at the initial time. But this means that $M_2(t)$ is a decreasing function with a negative time derivative bounded from above, so $M_2(t)$ hits $0$ at a time $t_0<\infty$. Due to nonnegativity of $u$ and mass conservation, this is in contradiction with global existence of a solution.

\qed 

\begin{rem}
Notice, that analyzing a proof of Theorem \ref{Tw.11} we conclude that once we pick up a proper $u_0$ then whatever is a $v_0$ we start from, anyway our solution must become unbounded. This is unlike in a one-dimensional setting of Theorem \ref{Tw.1}, where $v_0$ had to be chosen carefully so that (\ref{v_0}) is satisfied.  
\end{rem}

\vspace{0.3cm}
\noindent
\textbf{Acknowledgement.} T. Cie\'slak acknowledges a support from the Polish Ministry of Science and Higher Education under grant Iuventus Plus. The author is also grateful to Christian Stinner from Kaiserslautern for reading the early version of this paper and pointing out some misprints.


\begin{thebibliography}{99}
%
\bibitem{amann}
  \sc H.~Amann,
  \it Dynamic theory of quasilinear parabolic systems. III. Global existence.
  \rm Mathematische Zeitschrift {\bf 202}, 219-250 (1989).
\bibitem{BiN94}
	\sc P. Biler, T. Nadzieja, 
	\it Existence and nonexistence of solutions for a model of gravitational interaction of particles, I.
	\rm Colloq. Math. {\bf 66}, 319--334 (1994).
\bibitem{child} 
	\sc S.Childress, J.K. Percuss, 
	\it Nonlinear aspects of chemotaxis. 
	\rm Math. Biosci. {\bf 56}, 217-237 (1981).
\bibitem{CRAS}
	\sc T.~Cie\'slak, Ph.~Lauren\c cot,
	\it Finite-time blow-up of radially symmetric solutions to a critical quasilinear Smoluchowski-Poisson system. 
	\rm C.R. Math. Acad. Sci. Paris {\bf 347}, 237-242 (2009).
\bibitem{clKS} 
	\sc T.~Cie\'slak, Ph.~Lauren\c cot, 
	\it Finite time blow-up for a one-dimensional quasilinear parabolic-parabolic chemotaxis system.
	\rm Ann. Inst. H. Poincar\'e Anal. Non Lin\'eaire {\bf 27}, 437-446 (2010).	  
\bibitem{clSP} 
	\sc T. Cie\'slak, Ph. Lauren\c cot,  
	\it Looking for critical nonlinearity in the one-dimensional quasilinear Smoluchowski-Poisson system.
	\rm Discrete Contin. Dynam. Systems A {\bf 26}, 417-430 (2010).
\bibitem{clGlob_vs_blow} 
	\sc T. Cie\'slak, Ph. Lauren\c cot, 
	\it Global existence vs. blowup for the one dimensional quasilinear Smoluchowski-Poisson system. 
	\rm Progress in Nonlinear Differential Equations and Their Applications, {\bf vol. 80},  95-109 (2011).
\bibitem{ciesti}
	\sc T.~Cie\'slak, C.~Stinner,
	\it Finite-time blowup and global-in-time unbounded solutions to a parabolic-parabolic quasilinear Keller-Segel system in higher dimensions.
	\rm J.~Differential Equations {\bf 252}, 5832--5851 (2012).
\bibitem{ciesti2} 
	\sc T. Cie\'slak, C. Stinner
	\it Finite-time blowup in a supercritical quasilinear parabolic-parabolic Keller-Segel system in dimension $2$.
	\rm arXiV:1201.3270
\bibitem{CW08}
	\sc T. Cie\'slak, M. Winkler, 
	\it Finite-time blow-up in a quasilinear system of chemotaxis. 
	\rm Nonlinearity {\bf 21}, 1057--1076 (2008).
\bibitem{hil:volume}
	\sc T. Hillen, K. Painter, 
	\it Volume Filling and Quorum Sensing  in Models for Chemosensitive Movement. 
	\rm Canadian Applied Mathematics Quarterly {\bf 10}, 501-543 (2002).
\bibitem{HV:1}
	\sc M.A.~Herrero, J.J.L.~Vel\'{a}zquez, 
	\it A blow-up mechanism for a chemotaxis model. 
	\rm Ann. Scuola Norm. Sup. Pisa Cl. Sci. (4) {\bf 24}, 633-683 (1997). 
\bibitem{Ho02}
	\sc D. Horstmann, 
	\it On the existence of radially symmetric blow-up solutions for the Keller-Segel model.
	\rm J. Math. Biol. {\bf 44}, 463--478 (2002).
\bibitem{jl:expl}
	\sc W.~J\"ager, S.~Luckhaus, 
	\it On explosions of solutions to a system of partial differential equations modelling chemotaxis. 
	\rm Trans. Amer. Math. Soc. {\bf 329}, 819-824 (1992).	
\bibitem{KS:1}
	\sc E.F.~Keller, L.A.~Segel, 
	\it Initiation of slime mold aggregation viewed as an instability.
	\rm J. Theoretical Biology {\bf 26}, 399-415 (1970).
\bibitem{Na95}
	\sc T.~Nagai, 
	\it Blow-up of radially symmetric solutions to a chemotaxis system.
	\rm Adv. Math. Sci. Appl. {\bf 5}, 581-601 (1995).
\bibitem{Nagai}
	\sc T.~Nagai,
	\it Blowup of nonradial solutions to parabolic-elliptic systems modeling chemotaxis in two-dimensional domains.
	\rm J. Inequal. Appl. {\bf 6}, 37-55 (2001).
\bibitem{NSY97}
	\sc T. Nagai, T. Senba and K. Yoshida, 
	\it Application of the Trudinger-Moser inequality to a parabolic system 	of chemotaxis.
	\rm Funkcial. Ekvac. {\bf 40}, 411--433 (1997).
\bibitem{Nan}
	\sc V. Nanjundiah, 
	\it Chemotaxis, signal relaying and aggregation morphology.
	\rm J. Theoretical Biology  {\bf 42}, 63-105 (1973).
\bibitem{Nasri}
	\sc E. Nasreddine,
	\it Global existence of solutions to a parabolic-elliptic Keller-Segel system with critical degenerate diffusion.
	\rm arXiV 1207.4453
\bibitem{ss1}
	\sc T.~Senba, T.~Suzuki,
	\it Parabolic system of chemotaxis: blowup in a finite and the infinite time.
	\rm Methods Appl. Anal. {\bf 8}, 349-368 (2001).
\bibitem{win_mmas}
 \sc M.~Winkler, 
 \it Does a `volume-filling effect' always prevent chemotactic collapse?
 \rm Math.~Meth.~Appl.~Sci. {\bf 33}, 12-24 (2010).
\bibitem{win_bu}
 \sc M.~Winkler,
 \it Finite-time blow-up in the higher-dimensional parabolic-parabolic Keller-Segel system.
 \rm preprint, arXiv:1112.4156v1 (2011).                                        
\end{thebibliography}
\end{document}